\documentclass[conference]{IEEEtran}
\IEEEoverridecommandlockouts
\usepackage[left=0.63in, right=0.64in, top=0.7in, bottom = 0.99in]{geometry}
\usepackage{setspace}
\usepackage{cite}
\usepackage{subcaption}
\usepackage{amsmath,amssymb,amsfonts}
\usepackage{algorithm, algpseudocode}
\usepackage{graphicx}
\usepackage{subfig} 
\usepackage{textcomp}
\usepackage{xcolor}
\usepackage{tikz}
\usepackage{mathtools}
\usepackage{bm}
\usepackage{subfig}
\def\BibTeX{{\rm B\kern-.05em{\sc i\kern-.025em b}\kern-.08em
    T\kern-.1667em\lower.7ex\hbox{E}\kern-.125emX}}
\begin{document}

\title{\huge Movable Antennas in Wireless Systems: A Tool for Connectivity or a New Security Threat?}

\author{Youssef Maghrebi, Mohamed Elhattab, Chadi Assi, Ali Ghrayeb, Georges Kaddoum \vspace{-1cm}}

\maketitle
\allowdisplaybreaks
\begin{abstract}
The emergence of movable antenna (MA) technology has marked a significant advancement in the field of wireless communication research, paving the way for enhanced connectivity, improved signal quality, and adaptability across diverse environments. By allowing antennas to adjust positions dynamically within a finite area at transceivers, this technology enables more favourable channel conditions, optimizing performance across applications like mobile telecommunications and remote sensing. However, throughout history, the introduction of every new technology has presented opportunities for misuse by malicious individuals. Just as MAs can enhance connectivity, they may also be exploited for disruptive purposes such as jamming. In this paper, we examine the impact of an MA-enhanced jamming system equipped with $M$ antennas in a downlink multi-user communication scenario, where a base station (BS) with $N$ antennas transmits data to $K$ single-antenna users. Simulation results show that an adversary equipped with MAs reduce the system sum rate by 30\% more effectively than fixed-position antennas (FPAs). Additionally, MAs increase the outage probability by 25\% over FPAs, leading to a 20\% increase in the number of users experiencing outages. The highlighted risks posed by unauthorized use of this technology, underscore the urgent need for effective regulations and countermeasures to ensure its secure application.
\end{abstract}


\section{Introduction}
The advent of movable antenna technology represents an innovative shift in wireless communication systems, offering dynamic positioning capabilities to manage and improve signal propagation and reception adaptively. Unlike traditional fixed-position antennas (FPAs), movable antennas (MAs) can adjust their location within a designated region, exploiting spatial degrees of freedom (DoF) allowing them to optimize connectivity based on environmental and system-specific parameters. These systems harness spatial diversity more effectively than FPAs, particularly in scenarios with low mobility or limited diversity options, such as Industrial IoT, smart home applications and satellite communications \cite{survey-MA}. While the implementation of MAs may vary, they are generally defined as antennas connected to the radio frequency (RF) chain through a flexible cable, allowing mobility facilitated by a positioning mechanism or driver \cite{uplink-MA}. An early approach to MA design involved a reconfigurable linear array in which dipole antennas, mounted on stepper motors, could move along semicircular paths \cite{First-MA}. Another initial design allowed for the rotation of a uniform linear array (ULA) while maintaining fixed antenna spacing, aiming to achieve line-of-sight (LoS) transmission capacity \cite{second-MA}. The existing designs in the literature can be broadly classified into two categories: mechanical MAs and fluid antennas, also known as fluid antenna systems (FAS). Mechanical MAs utilize actuators, such as stepper motors, to enable antenna movement in three-dimensional space \cite{uplink-MA}. In contrast, fluid antennas are constructed using liquid metals or ionized solutions \cite{FAS}. These antennas are designed to move between discrete, predefined positions, referred to as ports, along a one-dimensional line. 

The literature on MAs primarily emphasizes their potential for enhancing link reliability, reducing interference, and supporting flexible deployment in challenging propagation environments. In \cite{field-response-model}, the authors compare the maximum channel gain of a single receive MA to that of its FPA counterpart, demonstrating that the MA system achieves significant performance improvements over conventional FPA systems, both with and without antenna selection (AS). The study in \cite{uplink-MA} showed that the total transmit power of multiple access systems can be significantly
decreased when using MAs compared to conventional FPA systems employing
AS. MAs can also significantly enhance physical layer security (PLS). In \cite{security1}, the authors demonstrate that an MA array considerably improves the secrecy rate compared to a conventional FPA array when transmitting confidential information in the presence of multiple eavesdroppers. Furthermore, \cite{security2} explores scenarios where perfect CSI about the eavesdroppers is unavailable, focusing on the transmission of confidential data from an MA-enabled array to multiple single-antenna receivers. The results reveal that MA-enabled PLS systems outperform FPAs in maximizing secrecy rates. 

While previous research highlights the significant advantages of MAs, this technology is also susceptible to misuse, as has been the case with many innovations throughout history. Hackers and malicious actors may attempt to exploit MAs for their harmful purposes. Therefore, it is crucial to examine the potential implications of MAs when used maliciously, a topic that is yet to be explored in the literature. One significant threat is jamming, where a malicious user aims to intensify interference against legitimate users to disrupt or stop their communication. In this paper, we investigate how MAs can enhance jamming efforts by a malicious user targeting legitimate users. Our main contributions are as follows:
\begin{itemize}
    \item We analyze the impact of an MA-assisted jammer in a downlink multiple-input single-output (MISO) system, where a base station (BS) serves $K$ users through space division multiple access (SDMA).
    \item We formulate an optimization problem in which the jammer decides both the antenna locations and beamforming vectors to minimize the overall system sum rate. Due to its non-convexity, the problem is divided into two sub-problems, which are alternately solved until convergence.
    \item We perform numerical examples under different scenarios, comparing the system sum rate and outage probability achieved with MAs versus FPAs. Results indicate that MAs reduce the system sum rate 30\% more than FPAs. Moreover, MAs raise the outage probability by 25\% compared to FPAs, resulting in a 20\% increase in the number of users affected by outages.
\end{itemize}
The remainder of this paper is structured as follows. Section II outlines the considered system model as well as the adopted channel models. In section III, we describe the problem formulation and detail the solution approach. Section IV is dedicated to the discussion and analysis of simulation results. Lastly, conclusions are drawn in Section V. In this paper, vectors and matrices are denoted by bold-face lower-case and upper-case letters, respectively. We denote \( \|\bm{M}\|_F \) as the Frobenius norm of any general matrix \( \bm{M} \). 

\section{System Model}

\subsection{Network Model}

\begin{figure}[t] 
  \centering
  \includegraphics[width=0.6\linewidth]{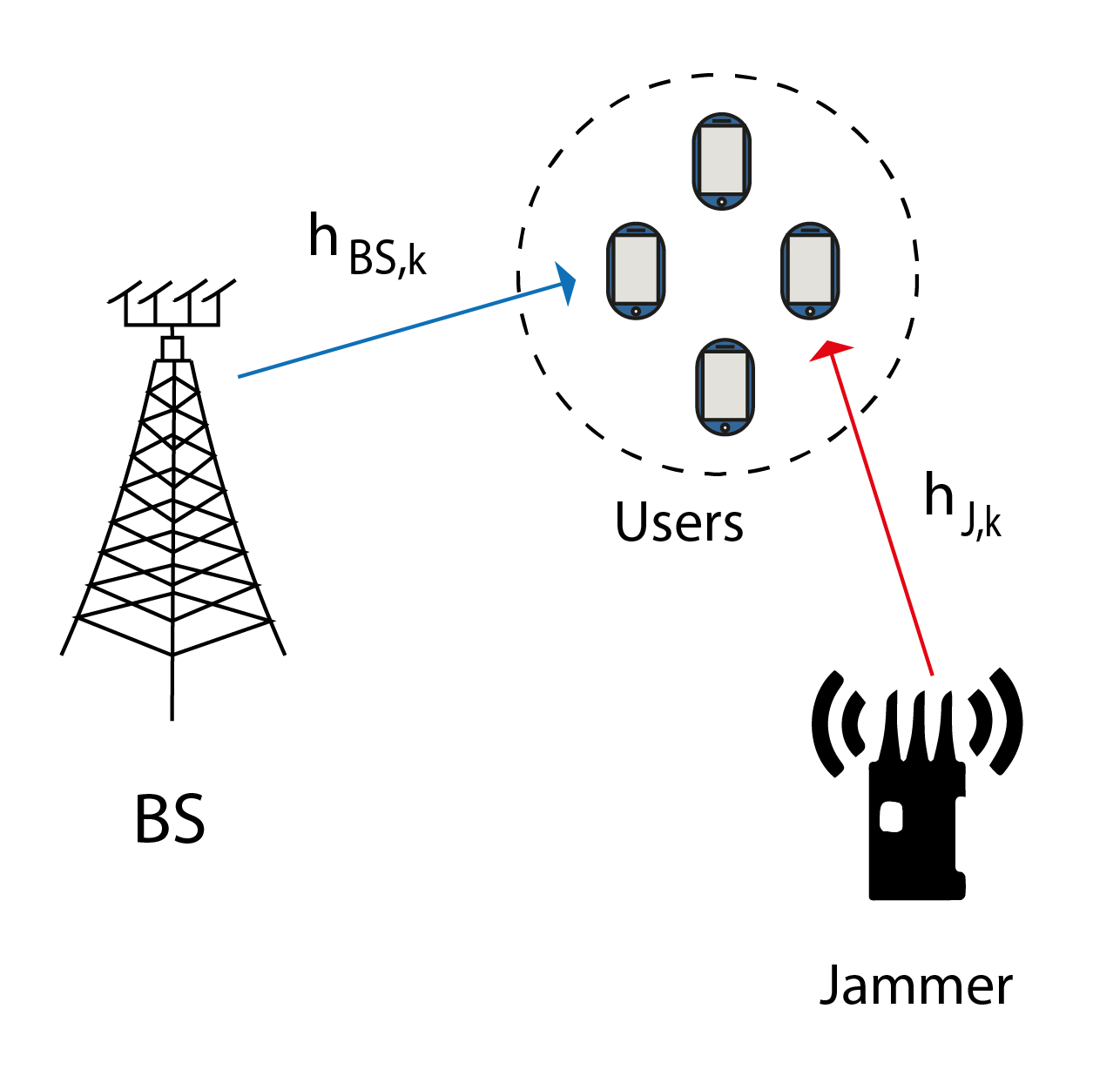} 
  \caption{System model}
  \label{fig:system model} 
\end{figure}

As shown in Fig.~1, a BS equipped with a ULA of size \(N\) serves \(K \leq N\) single-antenna users in downlink transmission at frequency \( f = \frac{c}{\lambda_{\text{BS}}} \). Meanwhile, a jammer attempts to disrupt this communication using a device with \(M\) movable antennas, as illustrated in Fig.~2. Though structurally similar to a ULA, the jammer’s antennas are mobile and can adjust positions according to the jammer's objectives. Each antenna has defined movement boundaries in the \((X, Z)\) plane, ranging from \([x_{\text{min}}, x_{\text{max}}] \times [z_{\text{min}}, z_{\text{max}}]\). Along the \(y\)-axis, the antennas are unrestricted provided that the spacing between any two consecutive antennas remains at least \(D = 2\lambda_{\text{J}}\) to avoid coupling effects~\cite{coupling-effect}, where \(\lambda_{\text{J}}\) denotes the jammer’s transmitted signal wavelength. The direct channel from the BS to user \(k\) is represented by \(\bm{h}_{\text{BS},k} \in \mathbb{C}^{N\times 1}\), and the channel from the jammer to user \(k\) is denoted by \(\bm{h}_{\text{J},k} \in \mathbb{C}^{M\times 1}\)\hspace{0.1cm}$\forall k \in \mathcal{K}\triangleq\{1, \ldots, K\}$. We denote $n_k \sim \mathcal{CN}(0, \sigma^2_k) \hspace{0.1cm}\forall k \in \mathcal{K}$, as the additive white Gaussian noise (AWGN).

Unaware of the jammer's presence, the BS optimizes its beamforming vectors based solely on the BS-user channels. Consequently, the BS's transmitted baseband signal can be expressed as $x = \sum_{k \in \mathcal{K}} \bm{w}_ks_k$, where $\bm{w}_k \in \mathbb{C}^{N\times 1}$ and $s_k$ are the BS beamforming vector and the information symbol for user $k$ with zero mean and unit variance. The jammer aims to disrupt the communication between the BS and the users by transmitting a signal \( z \) through its antenna array, expressed as \( z = \sum_{k \in \mathcal{K}} \bm{v}_k q_k \), where \( \bm{v}_k \in \mathbb{C}^{M \times 1} \) is the jammer's beamforming vector, and \( q_k \) is the jamming signal symbol directed toward user \( k \). We suppose that the jammer transmits over the same frequency used by the BS, i.e $\lambda_{BS} = \lambda_J = \lambda$. Additionally, the jammer possesses channel state information (CSI) for all channels between itself and the users\footnote{We assume the jammer performs prior channel estimation, storing estimated channels for various user locations in a dictionary. During jamming, it selects channels based on users' current locations. The effect of imperfect field-response information (FRI) due to minor user location misalignment is minimal, allowing the algorithm to maintain robust performance \cite{uplink-MA}.}; however, the jammer lacks the CSI for the BS-user links. Given the above, the actual received signal at user $k$ is expressed as\footnote{In this paper, we assume that \(M \geq K\), allowing the jammer to design a distinct beamforming vector for each user. If \(K > M\), the jammer can still effectively target a group of \(M\) users, making the problem equivalent to the scenario discussed in this paper.}:
\begin{equation}
y_k = \bm{h}_{BS,k}^H \bm{w}_k + \bm{h}^H_{J,k} \bm{v}_k + n_k    
\end{equation}

The rate expression for each user depends on the multiple access technique employed by the BS. In our case, the BS employs SDMA, thus the achievable rate for each user $k$ is given by: 
\begin{equation}
\small
\hspace{-0.1cm}R_k = \log_2 \left(1 + \frac{|\bm{h}_{\text{BS},k}^H \bm{w}_k|^2}{\sum_{\substack{j \in \mathcal{K} \\j \neq k}} |\bm{h}_{\text{BS},k}^H \bm{w}_j|^2 + \sum_{j\in \mathcal{K}} |\bm{h}_{\text{J},k}^H \bm{v}_j|^2 + \sigma_k^2} \right)
\end{equation}
where the term \(\sum_{\substack{j\in \mathcal{K}\\j \neq k}} |\bm{h}_{\text{BS},k}^H \bm{w}_j|^2\) represents the interference from other users, \( \sum_{\substack{j\in \mathcal{K}}} |\bm{h}_{\text{J},k}^H \bm{v}_j|^2 \) is the interference from the jammer, and \(\sigma_k^2\) is the noise power for user \(k\).

\begin{figure}[t] 
  \centering
  \includegraphics[width=\linewidth]{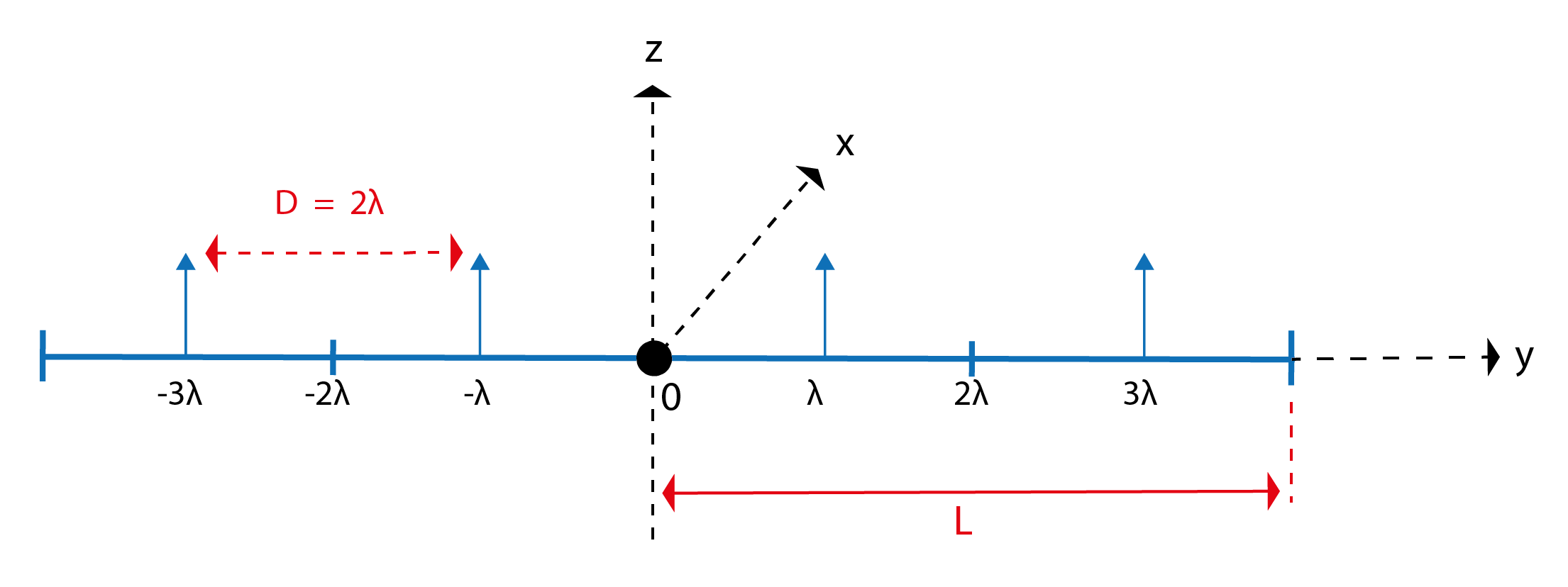} 
  \caption{Movable antenna array}
  \label{fig: array architecture} 
\end{figure}

\subsection{Channel Model}
In this paper, we consider narrow-band channels with slow fading. The BS-user links are modeled as Rayleigh fading, expressed as:
\begin{equation}
\bm{h}_{BS,k} = \sqrt{PL(d_{BS,k})}\bm{h}^{NLOS}_{BS,k} \hspace{0.1cm}\forall k \in \mathcal{K},   
\end{equation}
where $PL(d_{BS,k})$ is the distance dependent path loss modeled as $PL(d_{BS,k}) = \rho_o d^{-\alpha_{BS}}$. Here, $\rho_o$ represents the path loss at the reference distance of 1 m, $d_{BS,k}$ signifies the BS-user distance, and $\alpha_{BS}$ denotes the path loss exponent for the communication link. $\bm{h}^{NLOS}_{BS,k}$ represents the random non-LoS component modelled using a distribution with zero mean and unitary variance. 

The channel vector between the jammer and each user is determined by the propagation environment and the position of the MA at the jammer which is characterized by the field-response channel model \cite{field-response-model}. In this model, because the movement regions for the jammer's antennas are significantly smaller than the signal propagation distance, we assume the far-field condition holds. With this assumption, a plane-wave model can accurately represent the field response from the jammer's MA region to each user. Specifically, the angle of departure (AoD), angle of arrival (AoA), and the amplitude of the complex coefficient for each channel path between the base station and each user remain constant regardless of the MA's position within its movement region, while only the phase of each channel path changes with the MA’s position.

Let \( L^t_k \) and \( L^r_k \), for \( 1 \leq k \leq K \), denote the total number of transmit and receive channel paths from the jammer to user \( k \), respectively. The elevation and azimuth AoDs for the \( j \)-th transmit path from the jammer to user \( k \) are denoted by \( \theta^t_{k,j} \) and \( \phi^t_{k,j} \) for \( 1 \leq j \leq L^t_k \). Similarly, the elevation and azimuth AoAs for the \( i \)-th receive path to user \( k \) are denoted by \( \theta^r_{k,i} \) and \( \phi^r_{k,i} \) for \( 1 \leq i \leq L^r_k \). For simplicity, we define the virtual AoDs and AoAs as: 
\begin{align}
    \begin{aligned}
        \vartheta^t_{k,j} &= \cos \theta^t_{k,j} \cos \phi^t_{k,j}, \hspace{0.1cm}\varphi^t_{k,j} = \cos \theta^t_{k,j} \sin \phi^t_{k,j}\\
        \omega^t_{k,j} &= \sin \theta^t_{k,j}, \hspace{0.2cm} 1 \leq k \leq K, \hspace{0.2cm} 1 \leq j \leq L_k^t,\\
    \end{aligned} \label{eq:transmitter_angles} \\
    \begin{aligned}
       \vartheta^r_{k,i} &= \cos \theta^r_{k,i} \cos \phi^r_{k,i},  \hspace{0.1cm}\varphi^r_{k,i} = \cos \theta^r_{k,i} \sin \phi^r_{k,i}\\
       \omega^r_{k,i} &= \sin \theta^r_{k,i}, \hspace{0.2cm} 1 \leq k \leq K, \hspace{0.2cm} 1 \leq j \leq L_k^r.
    \end{aligned} \label{eq:receiver_angles}
\end{align}

In light of the above, the transmit and receive field-response vectors (FRVs) for the channel from the jammer to user $k$ are obtained as follows.
\begin{align}
\mathbf{f}_k(\mathbf{p}_m) &= \begin{bmatrix} e^{j \frac{2\pi}{\lambda} \rho^t_{k,1}(\mathbf{p}_m)}, e^{j \frac{2\pi}{\lambda} \rho^t_{k,2}(\mathbf{p}_m)}, \dots, e^{j \frac{2\pi}{\lambda} \rho^t_{k,L^t_k}(\mathbf{p}_m)} \end{bmatrix}^T, \\    
\mathbf{g}_k(\mathbf{u}_k) &= \begin{bmatrix} e^{j \frac{2\pi}{\lambda} \rho^r_{k,1}(\mathbf{u}_k)}, e^{j \frac{2\pi}{\lambda} \rho^r_{k,2}(\mathbf{u}_k)}, \dots, e^{j \frac{2\pi}{\lambda} \rho^r_{k,L^r_k}(\mathbf{u}_k)} \end{bmatrix}^T,
\end{align}
where, $\bm{u}_k \in \mathbb{R}^{3 \times 1}$ is the antenna position vector at user $k$ and $\bm{p}_m \in \mathbb{R}^{3 \times 1}$ is the m-th antenna position vector at jammer. $\rho^r_{k,j}(\mathbf{u}_k) = X_k \vartheta^r_{k,j} + Y_k \varphi^r_{k,j} + Z_k \omega^r_{k,j}$ for $1 \leq j \leq L^t_k$ represents the difference in signal propagation distance for the $j$-th receive channel path between the user antenna position $\mathbf{u}_k$ and the origin of the local coordinate system at user $k$. This indicates that the phase difference of the coefficient of the $j$-th receive channel path for user $k$ between user position $\mathbf{u}_k$ and $O_k$ is given by $\frac{2\pi}{\lambda} \rho^r_{k,j}(\mathbf{u}_k)$.

Similarly, the transmit FRV, $\mathbf{f}_k(\mathbf{p}_m)$, characterizes the phase differences in all $L^t_k$ transmit paths to user $k$, where $\rho^t_{k,i}(\mathbf{p}_m) = x_m \vartheta^t_{k,i} + y_m \varphi^t_{k,i} + z_m \omega^t_{k,i}$ for $1 \leq i \leq L^t_k$ represents the difference of the signal propagation distance for the $i$-th transmit channel path between the m-th jammer-antenna position $\mathbf{p}_m$ and the origin of the local coordinate system at jammer $O_0$.

The path-response matrix (PRM), $\Sigma_k \in \mathbb{C}^{L^t_k \times L^r_k}$, represents the response between all the transmit and receive channel paths from $O_0$ to $O_k$, $\forall k \in \mathcal{K}$. Specifically, the entry in the $i$-th row and $j$-th column of $\Sigma_k$ is the response coefficient between the $i$-th transmit path and the $j$-th receive path for user $k$. Subsequently, the channel vector from the jammer to user $k$ can be expressed as:
\begin{equation}
    \mathbf{h}_{J,k}(\mathbf{P}) = \mathbf{F}_k^H(\mathbf{P}) \Sigma_k \mathbf{g}_k, \quad \forall k \in \mathcal{K},
\end{equation}
where $\mathbf{F}_k = \begin{bmatrix} \mathbf{f}_k(\mathbf{p}_1), \mathbf{f}_k(\mathbf{p}_2), \dots, \mathbf{f}_k(\mathbf{p}_M) \end{bmatrix} \in \mathbb{C}^{L^t_k \times M}$ is the field-response matrix (FRM) at the jammer targeting user $k$, and $\bm{P} = [\bm{p}_1, \bm{p}_2, \ldots \bm{p}_M] \in \mathbb{R}^{3\times M}$ is the antenna position matrix at the jammer. Since the positions of the antenna at each user are fixed, the field receive vectors $g_k, \forall k \in \mathcal{K}$ are constant.

\section{Problem Formulation and Solution Approach}
\subsection{Problem Formulation}
In this paper, the aim of the jammer is to jointly optimize the MA location matrix $\bm{P}$, while designing the jamming beamforming matrix $\bm{V} = [\bm{v}_1, \bm{v}_2, \ldots, \bm{v}_K] \in \mathbb{C}^{M \times K}$ to minimize the network sum rate, while adhering to available power budget $P_{J}$ and the antenna array hardware constraints, which can be formulated as:
\begin{subequations}
\begin{align}
\mathcal{P}: \underset{\substack{\bm{V}, \bm{P}}}\min \hspace{0.2cm} & \sum_{k \in \mathcal{K}} R_k \\
\text{s.t.} \hspace{0.2cm} & \sum_{k \in \mathcal{K}}\|\bm{v}_k\|^2 \leq P_J,  \label{eq:power_budget}\\
& -\lambda \leq x_m \leq \lambda \hspace{0.2cm}\quad \forall m \in \{1, \ldots, M\}, \label{eq:x_constraint}\\
& -\lambda \leq z_m \leq \lambda \hspace{0.25cm}\quad \forall m \in \{1, \ldots, M\}, \label{eq:z_constraint}\\
& -L \leq y_m \leq L \hspace{0.2cm}\quad \forall m \in \{1, \ldots, M\}, \label{eq:y_constraint1}\\
& y_m-y_{m-1} \geq 2\lambda \quad \forall m \in \{2, \ldots, M\} \label{eq:y_constraint2}.
\end{align}
\end{subequations}
Constraint \eqref{eq:power_budget} represents the jammer's power budget limit. Constraints \eqref{eq:x_constraint} and \eqref{eq:z_constraint} define the permitted movement region for the MAs within the \((X, Z)\) plane. The typical size of the antenna moving region is in the order of several wavelengths \cite{field-response-model}, however the most common value in research is a distance of $\pm\lambda$ around the $x$ and $z$ axes. As previously mentioned, the MAs in this architecture are free to move along the entire length \(2L\) of the array on the \(y\)-axis, as specified in \eqref{eq:y_constraint1}, provided they maintain a minimum spacing of \(2\lambda\) between each other, as ensured by \eqref{eq:y_constraint2}. Since the jammer lacks information about the BS-user links or the BS beamforming vectors, its objective is to align its beamforming with each jammer-user channel to maximize interference. Therefore, the problem is formulated as follows. 
\begin{subequations}
\vspace{-0.3cm}
\begin{align}
\mathcal{P}: \underset{\substack{\bm{V}, \bm{P}}}\min \hspace{0.2cm} & -\sum_{k \in \mathcal{K}} |\bm{h}^H_{J,k}(\bm{P}). \bm{v}_k|^2 \\
\text{s.t.} \hspace{0.2cm} & \eqref{eq:power_budget}-\eqref{eq:y_constraint2}.
\end{align}
\end{subequations} 

The problem $\mathcal{P}$ is a non-convex optimization problem for two main reasons. First, the dependence on antenna location matrix  $\bm{P}$ in the jammer-user channel $\bm{h}_{J,k}^H$ is non linear, making the objective function non-convex. Additionally, constraints \eqref{eq:x_constraint}-\eqref{eq:y_constraint1} are obviously non-convex. 

\subsection{Solution Approach}
Given the non-convexity of problem $\mathcal{P}$ and the coupling between variables $\bm{P}$ and $\bm{V}$ in the objective function, we divide the problem into two sub-problems that will be solved iteratively until convergence. In each sub-problem, we focus on optimizing a set of variables while fixing the rest, which is usually referred to as block coordinate descent (BCD) or alternating optimization (AO). The detailed algorithm is described in Algorithm \ref{alg:sca_algo}. 

\begin{algorithm}
\small
\caption{Joint Optimization of $\mathbf{V}$ and $\mathbf{P}$}
\label{alg:sca_algo}
\begin{algorithmic}[1]
\State \textbf{Input:} $\mathbf{p}_m^{(0)} = [x_m^{(0)}, y_m^{(0)}, z_m^{(0)}]^T \hspace{0.1cm}\forall m \in \mathcal{M}$,  $\mathbf{v}_k^{(0)} \hspace{0.1cm}\forall k \in \mathcal{K}$, $\lambda$, $\epsilon$, $T_{\text{1,max}}$, $T_{\text{2,max}}$, $\mathbf{g}_k  \forall k \in \mathcal{K}$, $\vartheta^t_{k,i}, \varphi^t_{k,i}, \omega^t_{k,i}$ for $1 \leq i \leq L^t_k, \forall k \in \mathcal{K}$
\State \textbf{Output:} $\mathbf{p}^{opt}_m \hspace{0.1cm}\forall m \in \mathcal{M}$, $\mathbf{v}^{opt}_k \hspace{0.1cm}\forall k \in \mathcal{K}$
\State Set $i = 0$

\While{$\|\mathbf{V}^{(i+1)} - \mathbf{V}^{(i)}\|_F > \epsilon$ and $\|\mathbf{P}^{(i+1)} - \mathbf{P}^{(i)}\|_F > \epsilon$ and $i < T_{\text{1,max}}$}
    \State Set $t = 0$, $\bm{V}^{(t)} = \bm{V}^{(i)}$
    \While{$\|\mathbf{V}^{(t+1)} - \mathbf{V}^{(t)}\|_F > \epsilon$ and $t < T_{\text{2,max}}$}
        \State Compute $\nabla_{\mathbf{V}} \Psi^{(t)}$ using the current value of $\bm{P}$
        \State Solve $\mathcal{P}_1$ to obtain $\bm{V}^{(t+1)}$ 
        \State \small{Update $\mathbf{V^{(t)}} = \mathbf{V}^{(t+1)}$}
        \State $t = t + 1$  
    \EndWhile
    \State \small{Update $\mathbf{V^{(i+1)}} = \mathbf{V}^{(t)}$}
    \State Set $t = 0$, $\bm{P}^{(t)} = \bm{P}^{(i)}$
    \While{$\|\mathbf{P}^{(t+1)} - \mathbf{P}^{(t)}\|_F > \epsilon$ and $t < T_{\text{2,max}}$}
        \State Compute $\nabla_{\mathbf{P}} \Phi^{(t)}$ using the current value of $\bm{V}$
        \State Solve $\mathcal{P}_2$ to obtain $\bm{P}^{(t+1)}$ 
        \State \small{Update $\mathbf{P} = \mathbf{P}^{(t+1)}$}
        \State $t = t + 1$  
    \EndWhile
    \State \small{Update $\mathbf{P^{(i+1)}} = \mathbf{P}^{(t)}$}
    \State $i = i + 1$
\EndWhile
\end{algorithmic}
\end{algorithm} 

\subsubsection{Optimizing the beamforming matrix $\bm{V}$} 
In this subproblem, the antenna position matrix $\bm{P}$ is fixed. Thus, the optimization problem can be reduced to the following problem: 
\begin{subequations}
\begin{align}
\mathcal{P}_1: \underset{\substack{\bm{V}}}\min \hspace{0.2cm} & \Psi(\bm{V}) \triangleq -\sum_{k \in \mathcal{K}} |\bm{h}^H_{J,k}(\bm{P}). \bm{v}_k|^2 \\
\text{s.t.} \hspace{0.2cm} & \eqref{eq:power_budget}.
\end{align}
\end{subequations} 
The function $\Psi$ is concave, thus to have a convex objective, we apply a first-order Taylor expansion around a point $\Psi^{(t)}$ as:
\begin{equation}
\Psi = \Psi^{(t)} + \nabla_{\mathbf{V}} \Psi^{(t)} \cdot (\mathbf{V} - \mathbf{V}^{(t)}),    
\end{equation}
where, \(\nabla_{\mathbf{V}} \Psi^{(t)}\) is the Jacobian of \(\Psi\) with respect to \(\mathbf{V}\), evaluated at \(\mathbf{V}^{(t)} = \begin{bmatrix} \mathbf{v}_1^{(t)}, \mathbf{v}_2^{(t)}, \dots, \mathbf{v}_K^{(t)} \end{bmatrix}\). The new formulation for the problem is:
\begin{subequations}
\begin{align}
\mathcal{P}_1: \underset{\substack{\bm{V}}}\min \hspace{0.2cm} & \nabla_{\mathbf{V}} \Psi^{(t)} \cdot (\mathbf{V} - \mathbf{V}^{(t)}) \\
\text{s.t.} \hspace{0.2cm} & \eqref{eq:power_budget}.
\end{align}
\end{subequations} 
This convex problem is solved iteratively with CVX \cite{cvx}, beginning from a feasible starting point $\bm{V}^{(0)}$ and continuing until convergence, which is determined by reaching a specified tolerance $\epsilon$ or the maximum number of iterations $T_{2,max}$.

\subsubsection{Optimizing the antenna-position matrix $\bm{P}$}
In this subproblem, the beamforming matrix $\bm{V}$ is now fixed, and the new optimization problem can be expressed as:
\begin{subequations}
\begin{align}
\mathcal{P}_2: \underset{\substack{\bm{P}}}\min \hspace{0.2cm} & \Phi(\bm{P}) \triangleq -\sum_{k \in \mathcal{K}} |\bm{h}^H_{J,k}(\bm{P}). \bm{v}_k|^2 \\
\text{s.t.} \hspace{0.2cm} & \eqref{eq:x_constraint}-\eqref{eq:y_constraint2}.
\end{align}
\end{subequations} 
To alleviate the non-convexity from the objective function, we use a similar technique to the first sub-problem as follows.
\begin{equation}
\Phi = \Phi^{(t)} + \nabla_{\mathbf{P}} \Phi^{(t)} \cdot (\mathbf{P} - \mathbf{P}^{(t)}),    
\end{equation}
where, \(\nabla_{\mathbf{P}} \Phi^{(t)}\) is the Jacobian of \(\Phi\) with respect to \(\mathbf{P}\), evaluated at \(\mathbf{P}^{(t)} = \begin{bmatrix} \mathbf{p}_1^{(t)}, \mathbf{p}_2^{(t)}, \dots, \mathbf{p}_M^{(t)} \end{bmatrix}\).
Constraint \eqref{eq:y_constraint2} is convex, however, \eqref{eq:x_constraint}-\eqref{eq:y_constraint1} need to be transformed into the following convex constraints:
\begin{align}
& x_m \leq \lambda, \quad -x_m \leq \lambda  \label{eq:modified_x}\\
& z_m \leq \lambda, \hspace{0.1cm}\quad -z_m \leq \lambda  \label{eq:modified_z}\\
& y_m \leq L, \quad -y_m \leq L \label{eq:modified_y}
\end{align}
Using the above approximations, problem $\mathcal{P}_2$ is written as:
\begin{subequations}
\begin{align}
\mathcal{P}_2: \underset{\substack{\bm{P}}}\min \hspace{0.2cm} & \nabla_{\mathbf{P}} \Phi^{(t)} \cdot (\mathbf{P} - \mathbf{P}^{(t)}) \\
\text{s.t.} \hspace{0.2cm} & \eqref{eq:y_constraint2}, \eqref{eq:modified_x}-\eqref{eq:modified_y}.
\end{align}
\end{subequations} 
Similar to the first sub-problem, $\mathcal{P}_2$ is a convex problem that can be iteratively solved with CVX, starting from a feasible point $\bm{P}^{(0)}$ and continuing until convergence. 
\section{Simulation Results and Discussion}
In this section, we analyze the performance of the proposed MA-enhanced jamming solution through simulations, focusing on its effectiveness in maximizing interference against legitimate users. We present key performance metrics to evaluate the impact of MAs in various jamming scenarios, comparing these results to conventional FPAs.

\subsection{Simulation Setup}
The simulation setup is as follows: the BS is positioned at \((0, 0)\), while users are uniformly distributed within a circular area of radius $40$ m centered at \((50, 50)\). The jammer is located at \((100, 0)\). All coordinates are given in meters. We adopt the PRM matrix model used in \cite{uplink-MA}, where the PRM for each user \( k \), \( \Sigma_k = \operatorname{diag} \{ \sigma_1, \sigma_2, \dots, \sigma_L \} \), is a diagonal matrix where each diagonal element \( \sigma_i \) follows a circularly symmetric complex Gaussian (CSCG) distribution, \( \mathcal{CN}(0, \frac{c_k^2}{L}) \). Here, \( c_k^2 = g_0 d_k^{-\alpha_J} \) represents the expected channel power gain for user \( k \), with \( g_0 \) denoting the average channel power gain at a reference distance of 1 meter, and \( \alpha_J \) being the path-loss exponent. The elevation and azimuth AoDs and AoAs for the channel paths of each user are modeled as random variables with a joint probability density function (PDF) given by \( f_{\theta_k,\phi_k}(\theta_k, \phi_k) = \frac{\cos \theta_k}{2\pi} \), where \( \theta_k \in \left[ -\frac{\pi}{2}, \frac{\pi}{2} \right] \) and \( \phi_k \in \left[ -\frac{\pi}{2}, \frac{\pi}{2} \right] \).
Key simulation parameters are summarized in Table~\ref{tab: settings} and remain consistent across all simulations unless specified otherwise.
\begin{table}[!t]
  \centering
  \scriptsize
  \begin{tabular}{|c|c|c|c|}
    \hline
    \textbf{Parameter} & \textbf{Value} & \textbf{Parameter} & \textbf{Value} \\
    \hline
    $\alpha_{BS},\alpha_{J}$ & 2.8 & $\rho_o$ & -30 dB\\
    $P_{BS}$ & 40 dBm & $g_o$ & -40 dB\\
    $P_{J}$ & 30 dBm  & $\sigma^2$ & -80 dBm\\
    K,N,M & 4 & $R_{th,k}$ & 1 bps/Hz\\
    $L_t,L_r$ & 6 & $\lambda$ & 0.01 m\\
    \hline
  \end{tabular}
  \caption{Simulation Parameters}
  \label{tab: settings}
  \vspace{-0.2cm}
\end{table}

\begin{figure*}[t] 
    \centering
    \begin{subfigure}[b]{0.48\textwidth}
        \centering
        \includegraphics[width=\textwidth]{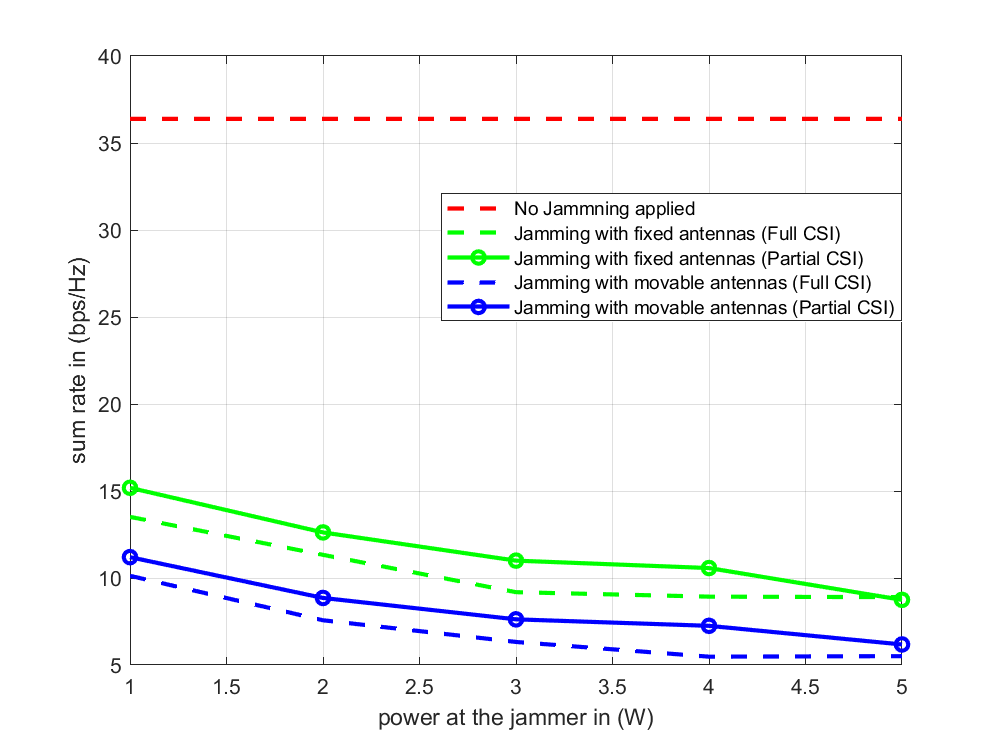}
        \caption{Sum rate change according to the available power at the jammer}
        \label{fig:actor_loss_1}
    \end{subfigure}
    \hfill
    \begin{subfigure}[b]{0.48\textwidth}
        \centering
        \includegraphics[width=\textwidth]{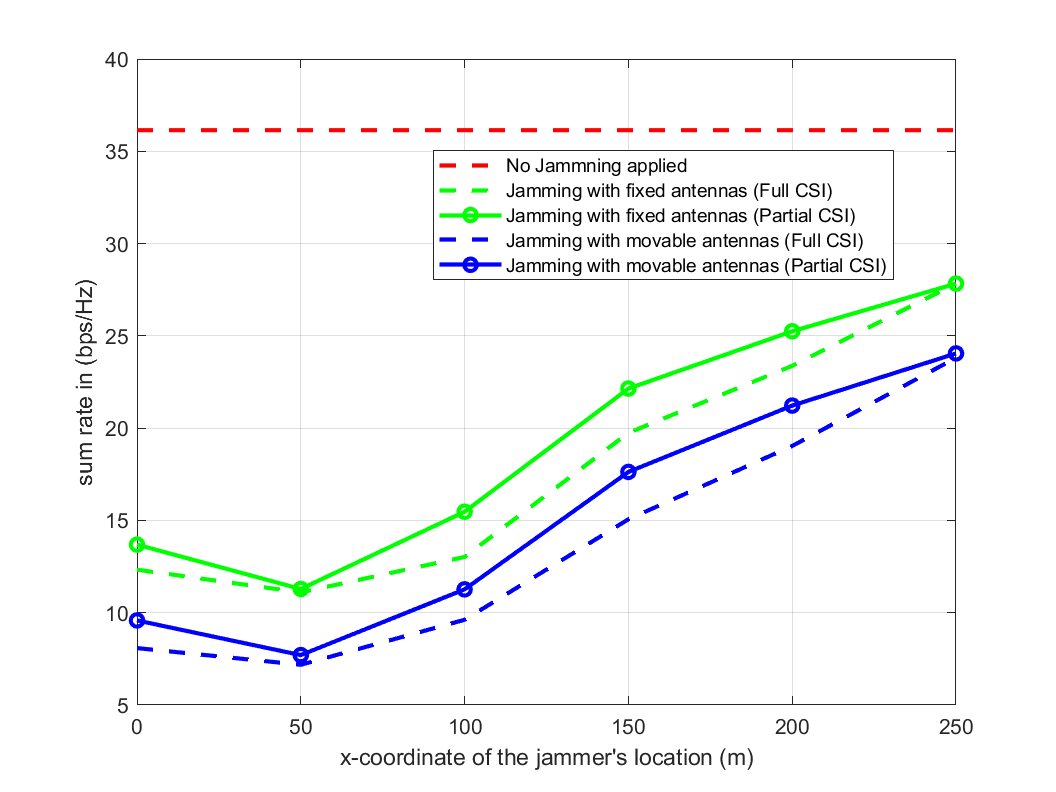}
        \caption{Sum rate change according to the jammer's location}
        \label{fig:critic_loss_1}
    \end{subfigure}
    \caption{Comparison of sum rate changes with respect to jammer power and location}
    \label{fig:combined_figures}
\end{figure*}

\begin{figure*}[t] 
    \centering
    \begin{subfigure}[b]{0.48\textwidth}
        \centering
        \includegraphics[width=\textwidth]{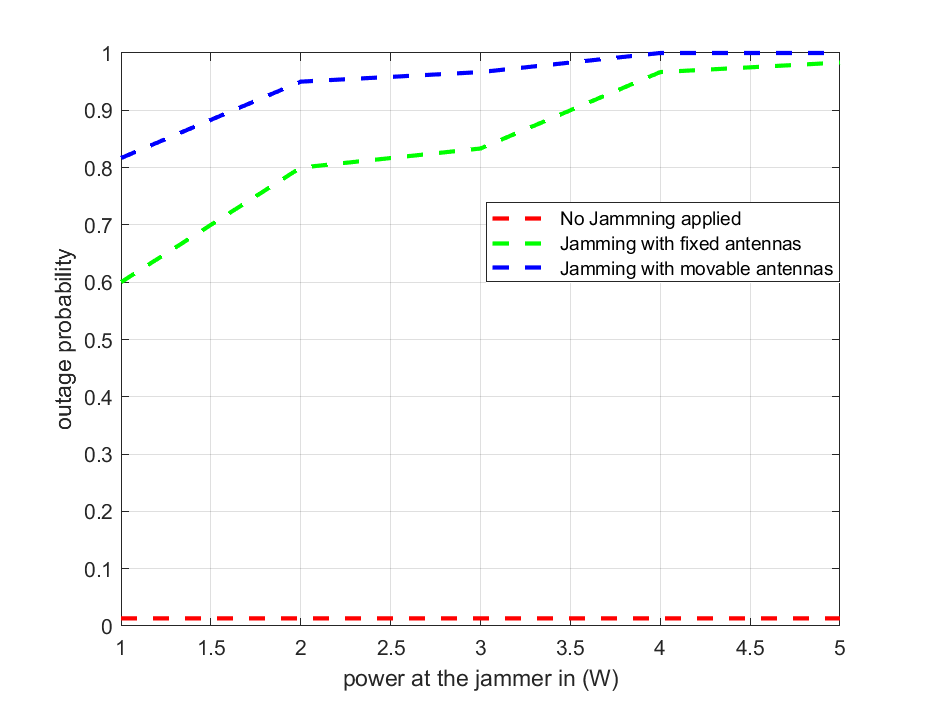}
        \caption{Outage probability versus the jamming power budget}
        \label{fig:actor_loss_1}
    \end{subfigure}
    \hfill
    \begin{subfigure}[b]{0.48\textwidth}
        \centering
        \includegraphics[width=\textwidth]{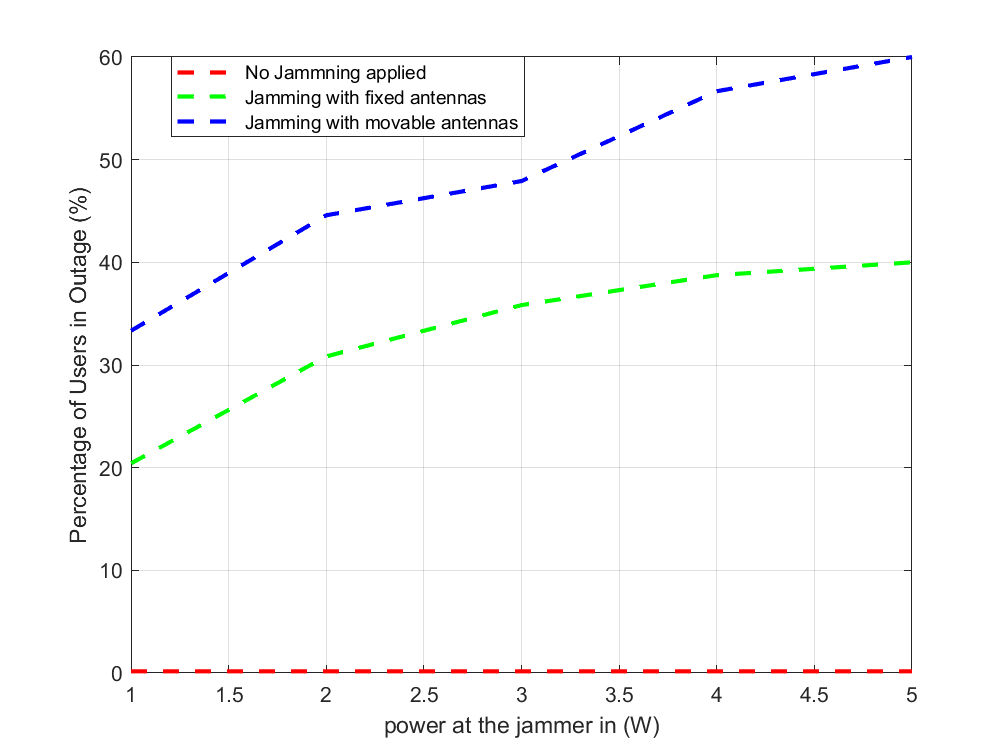}
        \caption{Percentage of users in outage versus the jamming power}
        \label{fig:critic_loss_1}
    \end{subfigure}
    \caption{Analysis of system outage probability with respect to jamming power}
    \label{fig:combined_figures}
\end{figure*}

\subsection{Results Analysis}
\subsubsection{Impact on the system sum-rate}
As previously mentioned, the jammer is assumed to possess only the CSI of its connections with each user. To evaluate the performance of jamming with MAs against conventional FPAs, while also considering the implications of lacking BS-user CSI and BS beamforming on the jammer's optimization, we define and compare the following baseline scenarios:

\begin{itemize} \item \textbf{Jamming with Fixed Antennas (Full CSI):} The jammer employs an FPA to disrupt communication while having access to both its jammer-user channels and complete BS information, which includes BS-user CSI and beamforming details. \item \textbf{Jamming with Fixed Antennas (Partial CSI):} In this scenario, the jammer utilizes an FPA but only has CSI about jammer-user channels. \item \textbf{Jamming with Movable Antennas (Full CSI):} Here, the jammer employs MA-enabled jamming, having full access to all system information. \item \textbf{Jamming with Movable Antennas (Partial CSI):} In this case, the jammer only possesses jammer-user CSI while using an MA array. \end{itemize}

The analysis of the impact of jamming on the sum rate performance, based on varying jammer locations and power is conducted in Fig~3. The simulations were conducted using the Monte-Carlo method, with results averaged over 100 different channel realizations. In the absence of jamming, the system achieves a stable sum rate of about 36 bps/Hz, serving as a baseline. When a jammer is introduced, its effectiveness varies significantly depending on its configuration.

In Fig.~3.a, we present the system sum rate achieved based on the power budget available at the jammer. As anticipated, the intensity of jamming, whether using MAs or FPAs, increases with the jammer's power. Notably, jamming with MAs proved to be, on average, $30\%$ more effective than using FPAs, resulting in an average reduction of 4 (bps/Hz) across all simulation points. Additionally, MA-based jamming results in a 15\% greater reduction in sum rate than its FPA counterpart compared to the baseline value when no jamming is applied. This significant difference could be a crucial factor in disrupting sensitive communication scenarios. Another noteworthy observation is that the difference in results between having partial CSI and full CSI is minimal, demonstrating the effectiveness of this jamming optimization despite its simplicity. In fact, the only requirement for the jammer is to estimate the CSI of its channels.   

The impact of the jammer's location on the system sum rate is illustrated in Fig.~3.b through the variation of the jammer's x-coordinate. The differences in values between jamming with MAs and FPAs continue to reflect the same gain percentages discussed in Fig.~3.a, thereby confirming the validity of the results. Jamming is most effective when the jammer is closest to the users, particularly near the center of the user distribution $(x_{jammer}=50~m)$. When positioned at equal distance from the center of the users, there is no significant difference in the achieved sum rate whether the jammer is located on the side of the BS $(x_{jammer}=0~m)$ or the opposite side $(x_{jammer}=100~m)$, indicating that the effectiveness of jamming primarily relies on the channel conditions between the jammer and the users. However, jamming from the BS side is slightly more effective, as the jammer can generate stronger interference with to BS-user links. This modest performance increase encourages the jammer to position itself on the opposite side of the BS to maximize its distance. Greater distance weakens the jammer's signal at the BS, making it harder to detect through signal strength-based techniques. Ultimately, as the jammer moves further from the center of the users, the intensity of jamming decreases.

\subsubsection{Impact on the system outage probability}
A critical factor in physical layer security is the outage probability $P_{\text{system outage}}$, defined as the probability of at least one user having a rate below the rate threshold. $P_{out,k}$ denotes the probability of outage of user $k$ given as:
\begin{equation}
P_{\text{out},k} = P(R_k < R_{\text{th}}), \hspace{0.2cm} \forall k \in \mathcal{K}.
\end{equation}
Since users are assumed to have independent channels, the probability that all users avoid outage (i.e., each user \( k \) meets \( R_k \geq R_{\text{th}} \)) is given by:
\begin{equation}
P_{\text{no outage}} = \prod_{k=1}^{K} \left( 1 - P_{\text{out},k} \right).
\end{equation}
Thus, the system outage probability is expressed as:
\begin{equation}
P_{\text{system outage}} = 1 - P_{\text{no outage}} = 1 - \prod_{k=1}^{K} \left( 1 - P_{\text{out},k} \right).    
\end{equation}
In Fig.~4, we evaluate the system outage probability when the jammer power varies between 1 and 5 W, through Monte Carlo simulation of 100 different channel realizations, where the rate threshold for each user is $R_{th} = 1~$bps/Hz. In the absence of jamming, the system maintains an outage probability of zero, as expected without interference. Fig.~4.a illustrates that jamming with fixed antennas yields an outage probability starting at 0.6 at 1 W and rising to approximately 0.95 by 5 W. However, jamming with MAs proves more disruptive, with the outage probability beginning around 0.75 at 1 W and approaching 1 by just 3.5 W. This demonstrates that movable antennas increase the outage probability more significantly—by about 0.2 across the power range—making MA arrays roughly 25\% more effective than FPAs at disrupting service for legitimate users.  

To deepen the analysis, we track the percentage of users experiencing outages among all users, as shown in Fig.~4.b. With jamming from FPAs, around 20\% of users are in outage at 1 W, rising moderately to 40\% at 5 W. In contrast, jamming with MAs starts at 30\% outage at 1 W and climbs sharply to about 60\% at 5 W, showing a linear-like increase in user disruptions. While the outage rate with FPAs begins to plateau near 40\%, MAs continue to increase steadily, with a performance gain that becomes more pronounced as power increases. This suggests that MAs cause more widespread outages, affecting roughly 15\% more users than FPAs at low power levels and up to 20\% more at 5 W. Overall, the higher outage rates with MAs underscore the importance of implementing effective countermeasures in network design.

\section{Conclusion}
In this paper, we have explored a primitive study on the potential misuse of movable antenna technology in wireless communication, highlighting its associated risks. Our investigation into an MA-enhanced jamming system revealed significant vulnerabilities that arise when such advanced technologies are exploited for malicious purposes. The results demonstrate that jamming using MAs is substantially more effective than FPAs, achieving, on average, a $30\%$ greater impact on system sum rate performance. Our analysis further illustrated that the jammer's location plays a crucial role in jamming efficacy, with proximity to users significantly enhancing the disruption potential. Notably, the difference in performance between scenarios with full channel state information (CSI) and partial CSI was minimal, indicating that even simple jamming strategies can effectively undermine communication integrity. Moreover, jamming with MAs significantly increases system outage probability compared to fixed antennas, with approximately 0.2 higher outage probability across a jammer power range of 1 to 5 W, making movable antennas about 25\% more effective at disrupting communication. Additionally, MAs lead to a greater proportion of users experiencing outages, affecting up to 20\% more users at higher power levels. These findings underscore the necessity of developing advanced security protocols that can adapt to the evolving landscape of wireless technology and counteract the threats posed by the unauthorized use of MAs in communication networks. 

Overall, while MA technology offers promising advancements in connectivity and performance, it is imperative that researchers and practitioners remain vigilant to the security challenges it introduces. Future work should focus on establishing comprehensive frameworks for protecting against the malicious exploitation of MAs, ensuring that the benefits of this technology can be harnessed without compromising the security and reliability of wireless communication systems.

\bibliographystyle{IEEEtran}
\bibliography{biblio}

\end{document}